\documentclass[12pt]{amsart}
\usepackage{epsfig}
\usepackage{amsfonts}
\usepackage{amssymb}
\usepackage{latexsym}

\newtheorem{theorem}{Theorem}[section]
\newtheorem{corollary}[theorem]{Corollary}
\newtheorem{lemma}[theorem]{Lemma}

\theoremstyle{definition}
\newtheorem{definition}[theorem]{Definition}
\newtheorem{remark}[theorem]{Remark}
\newtheorem{example}[theorem]{Example}

\begin{document}

\title{Combinatorics of the toric Hilbert scheme}
\author{Diane Maclagan}
\address{Department of Mathematics, University of California, Berkeley, Berkeley, CA 94720}
\email{maclagan@math.berkeley.edu}
\author{Rekha R. Thomas}
\address{Department of Mathematics, Texas A\&M University, 
College Station, TX 77843}
\email{rekha@math.tamu.edu}

\date{\today}
\maketitle

\begin{abstract}
The toric Hilbert scheme is a parameter space for all ideals with the
same multi-graded Hilbert function as a given toric ideal.
Unlike the classical Hilbert scheme, it is unknown whether toric
Hilbert schemes are connected. We construct a graph on
all the monomial ideals on the scheme, called the flip graph, and
prove that the toric Hilbert scheme is connected if and only if the
flip graph is connected. These graphs are used to exhibit curves in
$\mathbb P^4$ whose associated toric Hilbert schemes have arbitrary
dimension. We show that the flip graph maps into the
Baues graph of all triangulations of the point configuration defining
the toric ideal. Inspired by the recent discovery of a disconnected
Baues graph, we close with results that suggest the existence of a
disconnected flip graph and hence a disconnected toric Hilbert scheme.  
\end{abstract}

\section{Introduction}
Let $A = [a_1 \cdots a_n]$ be a $d \times n$ integer matrix of rank
$d$ such that $ker(A) \cap \mathbb N^n = \{ 0 \}$ and let ${\mathbb N
  A} := \{\sum_{i=1}^{n} m_ia_i \, : \, m_i \in {\mathbb N} \}
\subseteq {\mathbb Z^d}$ be the
non-negative integer span of the columns $a_1, \ldots, a_n$ of $A$.
The symbol ${\mathbb N}$ denotes the set of natural numbers including
zero. Consider the ${\mathbb Z^d}$-graded polynomial ring $S :=
k[x_1,\ldots,x_n]$ over a field $k$ with $\deg x_i := a_i$ for all $i$
and an ideal $I \subseteq S$ that is homogeneous with respect to
the grading by $\mathbb N A$, which we call {\em
  $A$-homogeneous}. The $k$-algebra $R=S/I$ is called an $A$-{\em
  graded algebra} if its Hilbert function $H_R(b):=\dim_k (R_b)$ is: 
$$H_R(b) = \left\{ \begin{array}{lll} 1 & \mathrm{ if }\, \,  b\in
    \mathbb N A & \\ 0 & \mathrm{ otherwise} & \\
\end{array} \right.$$
The presentation ideal $I$ is called an $A$-{\em graded ideal} and if 
$I$ is generated by monomials it is called a {\em monomial $A$-graded
  ideal}. 

$A$-graded algebras were introduced by Arnold \cite{Arnold} who
investigated matrices of the form $A = [1 \,\, p \,\, q]$ where $p$ and
$q$ are positive integers. A complete classification of all $A$-graded
algebras arising from one by three matrices can be found in
\cite{Arnold}, \cite{Korkina} and \cite{KPR}. The generalization
to $d$ by $n$ matrices is due to Sturmfels \cite{berndpreprint}. The 
canonical example of an $A$-graded ideal is the {\em toric ideal} of
$A$, denoted as $I_A$. {\em Initial ideals} of $I_A$ \cite{GBCP} are
also $A$-graded. 

In \cite{berndpreprint}, Sturmfels constructed a parameter space 
whose points are in bijection with the distinct $A$-graded ideals in
$S$. This variety is the underlying reduced scheme of the 
{\em toric Hilbert scheme} of $A$, denoted as $H_A$, which has been
defined recently by Peeva and Stillman \cite{PS2}, \cite{PS1}.
The classical Hilbert scheme parameterizes all homogeneous, 
saturated ideals in $S$ with a fixed Hilbert polynomial, where 
$S$ is graded by total degree. However, unlike classical Hilbert 
schemes which are known to be connected \cite{Hartshorne}, it is
unknown whether toric Hilbert schemes are connected. Several of the
techniques applied to classical Hilbert schemes cannot be used in the
toric situation. In particular, the multigraded Hilbert function used
to define $A$-graded ideals is not preserved under a change of
coordinates. See \cite{PS1} for further discussions. The only cases in
which $H_A$ is known to be connected are when $A$ has corank one
(i.e. $n-d=1$) or two. In the former case the connectivity is trivial,
and in the latter it follows from results in \cite{GP}.

In Section 2 we define a graph on all the monomial $A$-graded ideals
in $S$, called the {\em flip graph} of $A$, by defining an adjacency
relation among these ideals. This 
generalizes the notion of adjacency between two monomial initial
ideals of the toric ideal $I_A$, given by the edges of the {\em state
  polytope} of $I_A$ \cite{ST}. Our main result in Section 3 reduces the
connectivity of the toric Hilbert scheme to a combinatorial problem. 

\vspace{.1in}

\noindent{\bf Theorem 3.1.}
The toric Hilbert scheme $H_A$ is connected if and only if the flip
graph of $A$ is connected.  

\vspace{.1in}

The flip graph of $A$ provides information on the structure of
$H_A$. In Section 4 we use these graphs to prove that 
two by five matrices can have toric Hilbert schemes of arbitrarily high
dimension. The projective toric variety of such a matrix is a curve in
$\mathbb P^4$. 

\vspace{.1in}

\noindent{\bf Theorem 4.1.} For each $j \in {\mathbb N} \backslash
\{0\}$, there exists a two by five matrix $A(j)$ such that its toric
Hilbert scheme $H_{A(j)}$ has an irreducible component of dimension at
least $j$. 

\vspace{.1in}

In Section 5 we relate the flip graph of $A$ to the {\em Baues graph}
of $\mathcal A$ which is a graph on all the triangulations of the point
configuration ${\mathcal A} := \{ a_1, \ldots, a_n \} \subset {\mathbb
  Z^d}$ consisting of the columns of $A$. The edges of
the Baues graph are given by {\em bistellar flips}. This graph and its
relatives have been studied extensively in discrete geometry
\cite{reiner}. Sturmfels proved that the radical of a monomial
$A$-graded ideal $I$ is the {\em Stanley-Reisner} ideal of a
triangulation of $\mathcal A$, which we denote as
$\Delta(rad(I))$ (see Theorem~4.1 in \cite{berndpreprint} or 
Theorem~10.10 in \cite{GBCP}). This gives a map from the vertices  
of the flip graph into the vertices of the Baues graph. We extend this
map to the edges of the flip graph.

\vspace{.1in}

\noindent{\bf Theorem 5.2.} If $I$ and $I'$ are adjacent monomial 
$A$-graded ideals in the flip graph of $A$, then either they have the
same radical and hence $\Delta(rad(I)) = \Delta(rad(I'))$ or
$\Delta(rad(I))$ differs from $\Delta(rad(I'))$ by a bistellar flip.

\vspace{.1in}

Recently Santos \cite{Santos} constructed a configuration
${\mathcal A}$ with a disconnected Baues graph, settling 
the {\em generalized Baues problem} (see \cite{reiner} for a
survey). Although his example does not immediately give a disconnected
flip graph, it strongly supports the possibility of one. In Section 6
we explain this connection and provide results that point
toward a disconnected flip graph and hence, by Theorem 3.1, a
disconnected toric Hilbert scheme.    

\section{The Flip Graph of $A$}

In this section we define an adjacency relation on all monomial
$A$-graded ideals which, in turn, defines the flip graph of $A$. This 
graph is the main combinatorial object and tool in this paper. We first
recall the definition of an $A$-graded ideal.
 
\begin{definition}\label{agi}
  Let $A = [a_1 \cdots a_n] \in {\mathbb Z^{d \times n}}$ be a matrix
  of rank $d$ such that $ker(A) \cap \mathbb N^n = \{ 0 \}$ and let
  ${\mathbb N A} := \{ \sum_{i=1}^n m_ia_i \, : \, m_i \in {\mathbb N}
  \}$. An ideal $I$ in $S=k[x_1,\ldots,x_n]$ with $\deg x_i = a_i$ is
  called an $A$-graded ideal if $I$ is $A$-homogeneous and $R = S/I$
  has the ${\mathbb Z^d}$-graded Hilbert function: 
  $$H_R(b) := \dim_k (R_b) = \left\{ \begin{array}{lll} 1 & \mathrm
      { if }\, \,  b\in \mathbb N A & \\ 0 & \mathrm{ otherwise} & \\ 
\end{array} \right.
$$
\end{definition}

The canonical example of an $A$-graded ideal is the toric ideal $I_A$
which is the kernel of the ring homomorphism $\phi : S \rightarrow
k[t_1^{\pm}, \ldots, t_d^{\pm}]$ given by $x_j \mapsto t^{a_j}$. 
See \cite{GBCP} for more information. To see that $I_A$ is $A$-graded,
recall that $I_A = \langle x^u - x^v \, : \, Au = Av, \, u,v \in
{\mathbb N^n} \rangle$, and is hence $A$-homogeneous. For each $b \in
{\mathbb N A}$, any two monomials $x^u$ and $x^v$ in $S$ of $A$-degree
$b$ (i.e. with $Au = Av = b$) are $k$-linearly dependent modulo $I_A$
making $\dim_k((S/I_A)_b) = 1$. If $b \in {\mathbb Z^d} \backslash 
{\mathbb N A}$, $(I_A)_b$ is empty.

Given a {\em weight vector} $w \in {\mathbb N^n}$, the initial ideal 
of an ideal $I \subseteq S$ with respect to $w$ is the 
ideal $in_w(I) := \langle in_w(f) : f \in I \rangle$ where $in_w(f)$
is the sum of all terms in $f$ of maximal $w$-weight. 
Our assumption that $ker(A) \cap {\mathbb N}^n = \{0\}$ implies that there
is a strictly positive integer vector ${w'}$ in the row space of
$A$. Using the binomial description of $I_A$ given above, we then see 
that $I_A$ is homogeneous with respect to the grading $deg(x_i) =
w'_i$. Hence, the {\em Gr\"obner fan} of $I_A$ covers ${\mathbb R^n}$  
and each cell in this fan contains a non-zero non-negative integer
vector in its relative interior (see Proposition~1.12 in
\cite{GBCP}). Therefore, for any weight vector $w \in {\mathbb Z^n}$, 
the initial ideal $in_{w}(I_A)$ is well defined as it coincides 
with $in_{\bar w}(I_A)$ where $\bar w$ is a non-negative integer
vector in the relative interior of the Gr\"obner cone of $w$.
Since the Hilbert function is preserved when passing from an ideal to
one of its initial ideals, all initial ideals of $I_A$ are also
$A$-graded. 

If $M$ is a monomial $A$-graded ideal, then for each $b
\in {\mathbb N A}$ there is a unique monomial of degree $b$ that does
not lie in $M$ and is hence a {\em standard} monomial of
$M$. Definition~\ref{agi} implies that all $A$-graded ideals are 
generated by $A$-homogeneous {\em binomials} (polynomials with at most
two terms) since any two monomials of the same $A$-degree have to be
$k$-linearly dependent modulo an $A$-graded ideal.  

There is a natural action of the algebraic torus $(k^*)^n$ on $S$
given by $\lambda \cdot x_i=\lambda_ix_i$ for $\lambda \in (k^*)^n$. 

\begin{definition}
An $A$-graded ideal is said to be {\em coherent} if it is of the form
$\lambda \cdot in_w(I_A)$ for some $\lambda \in (k^*)^n$ and $w \in
{\mathbb Z}^n$. 
\end{definition}

We recall the definition of the {\em Graver basis} of $A$ \cite{GBCP}.
For $u,v \in {\mathbb N^n}$ we write $u < v$ if for each 
$i=1,\ldots,n$, $u_i \leq v_i$ and $u \neq v$. 

\begin{definition}
A binomial $x^u-x^v$ with $Au=Av$ is a {\em Graver binomial} 
if there do not exist $u^{\prime}, \, \, v^{\prime} \in \mathbb N^n$ with
$Au^{\prime}=Av^{\prime}$ and $u^{\prime} < u$, $v^{\prime}<v$.  The
collection of all Graver binomials is called the {\em Graver
basis}, $Gr_A$.
\end{definition}

The following lemma is a
strengthening of Lemma 10.5 in \cite{GBCP} and was also independently 
discovered by Peeva and Stillman (\cite[Proposition~2.2]{PS1}). 
The universal Gr\"obner basis of an ideal is the union of all the finitely many reduced Gr\"obner bases of the ideal.

\begin{lemma}\label{gensaregraver}
  Let $I$ be an $A$-graded ideal, and let $\mathcal
  G=\{x^{a_1}-c_1x^{b_1},\dots,\\ x^{a_k}-c_kx^{b_k} \}$ be the 
  universal Gr\"obner basis of $I$. Here the $c_i$
  may be zero and for each binomial, $x^{a_i}$ and $x^{b_i}$ are not
  both in $I$.  If $c_i=0$, choose $b_i$ so that $Aa_i=Ab_i$
  and $x^{b_i} \not \in I$.  Then for all $i$, $x^{a_i} - x^{b_i}$ is
  a Graver binomial. Hence, every minimal generator of
  $I$ is of this form.
\end{lemma}

\begin{proof}
  If $x^{a_i}-c_ix^{b_i} \in \mathcal G$, then there is some term order
  $\prec$ such that one of $x^{a_i}$ and $x^{b_i}$ is a minimal
  generator of $in_{\prec}(I)$, and the other is standard for
  $in_{\prec}(I)$.  Since $in_{\prec}(I)$ is also $A$-graded, it
  suffices to prove the lemma for monomial $A$-graded ideals, where
  $c_i=0$ for all $i$.
  
  Suppose there exist an $i$ such that $x^{a_i}-x^{b_i}$ is not a
  Graver binomial.  Then there exists $u,v \in \mathbb N^n$ with
  $Au=Av$ such that $u<a_i$ and $v < b_i$.  Since $I$ is $A$-graded,
  one of $x^u$ or $x^v$ is in $I$.  If we have $x^u \in I$ then
  $x^{a_i}$ would not be a minimal generator of $I$, and if $x^v \in I$
  then $x^{b_i}$ would not be standard. Therefore, $x^{a_i}-x^{b_i}$
  is a Graver binomial for all $i$.
\end{proof}

\begin{definition}
An $A$-homogeneous ideal $I$ in $S$ is {\em weakly $A$-graded} if
$H_{S/I}(b) \in \{0,1 \}$ for all $b \in \mathbb Z^d$, and
$H_{S/I}(b)=1 \Rightarrow b \in \mathbb N A$.
\end{definition}

\begin{lemma} \label{UGB}
Let $I$ be an ideal which contains a binomial of the form
$x^a-cx^b$ for every Graver binomial $x^a-x^b$.  Then $I$ is weakly
$A$-graded. 
\end{lemma}

\begin{proof}
  It suffices to prove that $M=in_{\prec}(I)$ is weakly $A$-graded,
  where $\prec$ is any term order, since $in_{\prec}(I)$ has the same 
  Hilbert series as $I$. If $x^a-x^b$ is a Graver binomial,
  then since there is some $c$ with $x^a-cx^b \in I$, one of $x^a$ and
  $x^b$ lies in $M$.  Let $x^u$ and $x^v$ be two monomials of degree
  $b$, and let $x^a-x^b$ be a Graver binomial with $x^a |x^u$ and $x^b
  |x^u$.  Since one of $x^a$ and $x^b$ lies in $M$, one of $x^u$ and
  $x^v$ lies in $M$.  It thus follows that there is at most one
  standard monomial of $M$ in each degree $b$, and so $M$ is weakly
  $A$-graded.
\end{proof}

We now define a ``flipping'' procedure on a monomial $A$-graded ideal
which transforms this ideal into an ``adjacent''
monomial $A$-graded ideal. The idea is motivated by a similar
procedure for toric initial ideals which we describe briefly. 

The distinct monomial initial ideals of $I_A$ are in bijection 
with the vertices of the state polytope of $I_A$, an
$(n-d)$-dimensional polytope in ${\mathbb R}^n$ \cite{ST}. Two initial
ideals are said to be adjacent if they are indexed by adjacent
vertices of the state polytope. The edges of the state polytope are
labeled by the binomials in the universal Gr\"obner basis of $I_A$,
$UGB_A \subseteq Gr_A$. 

Suppose $I$ and $I'$ are two adjacent monomial initial ideals of $I_A$
connected by the edge $x^a-x^b$. The closure of the outer normal cone
at the vertex $I$ (respectively $I'$) is the {\em Gr\"obner cone} $K$
(respectively $K'$) of $I$ (respectively $I'$), the interior of which
contains all the weight vectors $w$ such that $in_w(I_A) = I$
(respectively $in_w(I_A) = I'$). The linear span of the common facet
of $K$ and $K'$ is the hyperplane $\{u \in {\mathbb R^n} \, :\, (a-b)
\cdot u = 0 \}$. When $w$ is in the interior of $K$, $in_w(x^a-x^b) =
x^a$, $x^a$ is a minimal generator of $I$ and $x^b \not \in I$, and
when $w$ is in the interior of $K'$, $in_w(x^a-x^b) = x^b$, $x^b$ is a
minimal generator of $I'$ and $x^a \not \in I'$. For a $w$ in the
relative interior of the common facet of $K$ and $K'$, $in_w(x^a-x^b)
= x^a-x^b$. Hence passing from $I$ to $I'$ involves ``flipping'' the
orientation of the binomial $x^a-x^b$. No other binomial in $UGB_A$
changes orientation during this passage. See \cite{HuT} for
details. We extend this notion of ``flip'' to all monomial $A$-graded
ideals.

\begin{definition}
Let $I$ be a monomial $A$-graded ideal and $x^a-x^b$ a Graver binomial
with $x^a$ a minimal generator of $I$ and $x^b \not \in I$. 
We define $I_{flip}$, the result of flipping over this binomial, to be
$$I_{flip} := \langle x^c | \exists \, \, d : x^c-x^d \in Gr_A, x^c
\in I, x^d \not \in I, c \neq a \rangle + \langle x^b \rangle.$$
\end{definition}

\begin{lemma} \label{weaklyA}
The ideal $I_{flip}$ is weakly $A$-graded.
\end{lemma}

\begin{proof}
Let $x^{\alpha}-x^{\beta}$ be a Graver binomial, with $x^{\alpha}
\in I$.  By Lemma \ref{UGB} it suffices to show that either
$x^{\alpha} \in I_{flip}$, or $x^{\beta} \in I_{flip}$.  Since
$x^{\alpha} \in I$, there is some (possibly identical) Graver binomial
$x^{\alpha^{\prime}}-x^{\beta^{\prime}}$ with $x^{\alpha^{\prime}}$ a
minimal generator of $I$, and $x^{\beta^{\prime}} \not \in I$, and
$x^{\alpha^{\prime}} |x^{\alpha}$.  If $\alpha^{\prime} \neq a$, then
$x^{\alpha^{\prime}} \in I_{flip}$, and so $x^{\alpha}\in I_{flip}$.  If
$\alpha^{\prime}=a$, then $\beta^{\prime}=b$, so $x^{\beta} \in
I_{flip}$.
\end{proof}

As defined above, to construct $I_{flip}$ requires knowledge of the
entire Graver basis.  However, the local change algorithm in
\cite{HuT} can be used to construct $I_{flip}$.   

\begin{lemma} \label{localflip}
  The ideal $I_{flip}$ is the initial ideal with respect to $x^a \prec
  x^b$ of $W_{a-b}=\langle x^c | c \neq a, x^c$ is a minimal generator
  of $I \rangle + \langle x^b-x^a \rangle$.
\end{lemma}

We note first that this initial ideal is well-defined. The only
non-trivial S-pairs formed during its construction are those of a
monomial with $x^b-x^a$, in which case the result is a monomial
multiple of $x^a$, so there is never any question of what the leading
term of a polynomial is.  This means that $I_{flip}$ is in fact the
initial ideal of $W_{a-b}$ with respect to {\em any} term order in
which $x^a \prec x^b$. We call $W_{a-b}$ a {\em wall ideal}
since in the coherent situation, it is the initial ideal
of any weight vector in the relative interior of the common facet/wall
between the Gr\"obner cones of $I$ and $I_{flip}$ \cite{HuT}.

\begin{proof}
Let $K$ be the initial ideal of $W_{a-b}$ with respect to $x^a \prec
x^b$. We first show the containment $K \subseteq I_{flip}$.  Let $x^c$
be a minimal generator of $K$. If $x^c=x^b$, or $x^c$ is a minimal
generator of $I$ other than $x^a$, then $x^c \in I_{flip}$. 
So we need only consider the case that $c=ra+g$, where $r > 0$ and
$a,b \not \leq g$, as this is the only other form minimal generators
of $K$ can have. In order to show that $x^c$ is in $I_{flip}$, it 
suffices to show that $x^c-x^d$ is a Graver
binomial, where $x^d$ is the unique standard monomial of $I$ of the 
same $A$-degree as $x^c$. 

Suppose $x^c-x^d$ is not a Graver
binomial, so we can write $c=\sum u_i + g^{\prime}$, $d=\sum v_i +
g^{\prime}$, where for each $i$, $x^{u_i}-x^{v_i}$ is a Graver
binomial. Since $x^d \not \in I$, we must have $x^{u_i} \in I$ and
$x^{v_i} \not \in I$ for all $i$.  If $u_i \neq a$ for some $i$, this
would mean that $x^{u_i}$, and hence $x^c$, was in $I_{flip}$. We can
thus reduce to the case where $g'=g$ and $v_i = b$ for all $i$, and so 
$d=rb+g$. Now since $x^c$ is a minimal generator of $K$, there must be
some minimal generator, $x^{\alpha}$, of $I$ for which the result of the
reduction of the S-pair of $x^{\alpha}$ and $x^b-x^a$ is
$x^c$. The only binomial that can be used in the reduction is
$x^b-x^a$, and hence there exists $l,m \geq 0$ such that $l+m=r$ and 
$x^{la+mb+g}$ is the least common multiple of $x^b$ and $x^{\alpha}$.  
If $l \neq 0$, then $x^a | x^{la+mb+g} | x^{\alpha+b}$. Since $x^a$ and
$x^b$ have no common variables, we get that $x^a | x^{\alpha}$ which 
contradicts $x^{\alpha}$ being a minimal generator of $I$.
So we must have $l=0$ and $x^{rb+g}$ is the least
common multiple of $x^b$ and $x^{\alpha}$.  But this implies that 
$x^{rb+g}=x^d$ is a multiple of $x^{\alpha}$ and hence in $I$, which
is a contradiction.  Therefore, this case cannot occur and 
we conclude that $K \subseteq I_{flip}$. 

We now show the reverse inclusion.  Suppose $x^c$ is a minimal
generator of $I_{flip}$ not equal to $x^b$, and $x^c-x^d$ is the
corresponding Graver binomial with $x^d \not \in I$.  We may
assume that $x^c$ is a multiple 
of $x^a$, as otherwise it is a generator of $W_{a-b}$, and thus in $K$
automatically.  Write $c=ra+\gamma$, where $a \not \leq \gamma$.
Suppose that $x^{rb+\gamma} \not \in I$. Then $d=rb+\gamma$, so we
must have $\gamma=0$ and $r=1$ to preserve $x^c-x^d$ being a Graver
binomial.  But then $c=a$, contradicting $x^c$ being a minimal
generator of $I_{flip}$.  Thus $x^{rb+\gamma} \in I$, and so there is
some $\alpha \neq a$ with $x^{\alpha}$ a minimal generator of $I$ such
that $\alpha \leq rb+\gamma$.  This means that $x^{rb+\gamma} \in
W_{a-b}$, and so $x^{ra+\gamma}=x^c\in W_{a-b}$ because $x^b-x^a \in
W_{a-b}$. Any monomial in $W_{a-b}$ is in $K$, so we conclude that
$x^c \in K$. 
\end{proof}

\begin{definition}
We say that a binomial $x^{a} -x^{b}$ in the Graver basis is 
{\em flippable} for a monomial $A$-graded ideal $I$ if $x^a$ is a
minimal generator of $I$, $x^b \not \in I$ and the ideal
$I_{flip}$ obtained by flipping $I$ over $x^a-x^b$ is again a monomial
$A$-graded ideal.
\end{definition}

We now give a characterization of when a binomial is flippable.

\begin{theorem} \label{flipcrit}
  Let $I$ be a monomial $A$-graded ideal, and $x^a-x^b$ a Graver
  binomial. Then $x^a-x^b$ is flippable for $I$ if and only if $I$ is
  the initial ideal with respect to $x^b \prec x^a$ of the wall ideal
  $W_{a-b}=\langle x^c | c \neq a, x^c$ is a minimal generator of $I
  \rangle + \langle x^a-x^b \rangle$.
\end{theorem}

\begin{proof}
Since $W_{a-b}$ is $A$-homogeneous, $I$ is the initial ideal of
$W_{a-b}$ if and only if $W_{a-b}$ is an $A$-graded ideal.  But by
Lemma \ref{localflip} $I_{flip}$ is an initial ideal of $W_{a-b}$, so
is $A$-graded exactly when $W_{a-b}$ is.
\end{proof}

\begin{definition} The {\em flip graph} of $A$ has as its vertices 
all the monomial $A$-graded ideals in $S$. There is an edge labeled by
the Graver binomial $x^a-x^b$ between two vertices $I$ and $I'$, if
$I'$ can be obtained from $I$ by flipping over $x^a-x^b$.
\end{definition}

\begin{remark}
The edge graph of the state polytope of $I_A$ is a 
subgraph of the flip graph of $A$. Since the state polytope of 
$I_A$ is $(n-d)$-dimensional, this subgraph is $(n-d)$-connected 
and so every vertex in this subgraph has valency at least $n-d$.
\end{remark}

Let $Flips_A$ denote the set of binomials labeling the edges of the
flip graph of $A$. Since the edges of the state 
polytope of $I_A$ are labeled by the elements in $UGB_A$, 
we have $UGB_A \subseteq Flips_A \subseteq Gr_A$. 

\begin{remark} (i) Gasharov and Peeva \cite{GP} proved that all
monomial $A$-graded ideals of corank two matrices are coherent.
Hence, in this case, the flip graph of $A$ is precisely the edge graph
of the state polytope of $I_A$, which is a polygon since $n-d=2$, and
$UGB_A = Flips_A$. However, 
even in this case, $Flips_A$ may be properly contained in $Gr_A$: for
$A = [1\,\,3\,\,7]$, $UGB_A = Flips_A = \{a^2c-b^3, a^3-b, ac^2-b^5, 
b^7-c^3, c-a^7, ab^2-c \}$ while $Gr_A = Flips_A \cup
\{a^4b-c\}$.\\
(ii) For $A = [1 \,\,3\,\,4]$, $UGB_A = Flips_A = Gr_A = 
\{ac^2-b^3, a^2c-b^2, b^4-c^3, b-a^3, ab-c, a^4-c\}$.\\
(iii) For $A = [3\,\,4\,\,5\,\,13\,\,14]$, $UGB_A \subsetneq Flips_A
\subsetneq Gr_A$. 
In this case, $Flips_A \backslash UGB_A = \{
a^2bcd-e^2\}$ while $Gr_A \backslash Flips_A = \{d^4-bc^4e^2,
ad^3-bc^2e^2, e^3-b^6cd, b^3cd^3-e^4, e^3- a^2c^2d^2, e^2-ab^5c, 
e^3-ab^2cd^2, e^3-a^4bd^2 \}$.
\end{remark}

For fixed $A$, let $S_A$ be the intersection of all the 
monomial $A$-graded ideals in $S$ and let $P_A := \langle x^ax^b :
x^a-x^b \in Gr_A \rangle$. Then $P_A$ is contained, sometimes
strictly, in $S_A$ since for each Graver binomial $x^a-x^b$, at least
one of $x^a$ or $x^b$ belongs to each monomial $A$-graded ideal. 

\begin{lemma}
If $x^a-x^b \in Gr_A$ has at least one of 
$x^a$ or $x^b$ in $P_A$, then $x^a-x^b \in Gr_A \backslash
Flips_A$. The converse is false. 
\end{lemma}

\begin{proof} 
Suppose $x^a - x^b$ is a flippable binomial for a monomial $A$-graded
ideal $M$ such 
that $x^a \in M$ and $x^b \not \in M$. If $x^a \in P_A \subseteq
S_A$ then $x^a \in M_{flip}$ and if $x^b \in P_A \subseteq S_A$
then $x^b \in M$ both of which are contradictions. 
To see that the converse is false, consider $$A = \left(
\begin{array}{cccccc} 
2 & 1 & 0 & 1 & 0 & 0\\
0 & 1 & 2 & 0 & 1 & 0 \\
0 & 0 & 0 & 1 & 1 & 2 \end{array} \right)$$ which has 29 monomial
$A$-graded ideals all of which are coherent. The binomial
$x_1x_4x_6-x_2x_3x_5 \in Gr_A \backslash Flips_A$, but neither
$x_1x_4x_6$ nor $x_2x_3x_5$ 
lies in $P_A = \langle
x_1x_2^2x_4, x_1x_3^2x_6, x_4x_5^2x_6, x_1x_2x_3x_5, x_2x_3x_4x_5, 
x_2x_3x_5x_6, \\ 
x_1x_2^2x_5^2x_6, x_2^2x_3^2x_4x_6, x_1x_3^2x_4x_5^2, 
x_1x_2x_3x_4x_5x_6  
\rangle$. 
\end{proof}

\section{Connection to the toric Hilbert Scheme}
In this section we explain the relevance of flips for the toric
Hilbert scheme $H_A$. We begin by describing the toric Hilbert scheme.

A parameter space for the set of $A$-graded ideals was first described
by Sturmfels \cite{berndpreprint}.  Peeva and Stillman improved on
this construction by producing the toric Hilbert scheme of $A$
\cite{PS2}, \cite{PS1}, which they show satisfies an important
universal property.  It is a version of their equations we explain
below.

A degree $b \in \mathbb N A$ is a Graver degree if there is some
Graver binomial $x^{\alpha}-x^{\beta}$ with $A\alpha=A\beta=b$.
We denote by $b_1,\dots,b_N$ the Graver degrees and by $m_i$ the
number of monomials of degree $b_i$.  Let
$$X=\mathbb P^{m_1-1} \times \mathbb P^{m_2-1} \times \dots \times
\mathbb P^{m_N-1}.$$

We now describe $H_A$ as a subscheme of $X$.  The coordinates of each
$\mathbb P^{m_i-1}$ can be labeled by the monomials of degree $b_i$ as
$\{\xi_u : Au=b_i \}$. A point $p\in X$ corresponds to a weakly
$A$-graded ideal $I_p$ by the following procedure: For each pair $x^u, x^v$
of degree $b_i$, we place the binomial $\xi_v x^u-\xi_u x^v$ in $I_p$.
For each Graver binomial $x^{\alpha}-x^{\beta}$ there thus is a
binomial of the form $x^{\alpha}-cx^{\beta}$ in the resulting ideal,
where $c$ may be zero.  This is immediate except in the case that
$\xi_u=\xi_v=0$.  In that case, choose $w$ with $Aw=A\alpha$ such that
$\xi_w \neq 0$.  Then the binomial $\xi_wx^u-\xi_ux^w \in I_p$, so
$x^u \in I_p$, and so $x^u-0\cdot x^v$ is the required binomial.  
Lemma \ref{UGB} now implies that $I_p$ is weakly $A$-graded.

We note that the toric ideal $I_A$ corresponds to the point in $X$
with $\xi_u=1$ for all $u$.  A monomial $A$-graded ideal corresponds
to a point in $X$ where for each $1 \leq i \leq N$ there is exactly one
$\xi_{u_i}=1$, and $\xi_v=0$ if $v\neq u_i$ for some $i$.  In general,
if $Au=b_i$, then $x^u \in I_p$ exactly if $\xi_u=0$.

We now give equations for $H_A$, which
guarantee that the resulting ideals $I_p$ are in fact $A$-graded.  Let $B
\subset \mathbb NA$ be a finite collection of degrees such that if a
weakly $A$-graded ideal generated in Graver degrees is $A$-graded in
every degree in $B$, then it is $A$-graded.  We know that such a $B$
exists because of bounds given by Sturmfels \cite{berndpreprint} and
Peeva and Stillman \cite{PS2}.

For each $b \in B$ we construct the matrix $M_b$ whose $d_b$ rows are
labeled by the monomials of degree $b$.  The $n_b$ columns of $M_b$
are labeled by pairs $x^u, x^v$ of degree $b$ such that there is some
Graver binomial $x^{\alpha}-x^{\beta}$ such that
$u=v-\alpha+\beta$.  The corresponding column consists of
$\xi_{\alpha}$ in the $x^u$ row, $-\xi_{\beta}$ in the $x^v$ row, and
zeroes elsewhere.

The global equations for $H_A$ are now given by the maximal minors of
$M_b$ for every $b\in B$.  To see that these equations guarantee that
$I_p$ is $A$-graded, note that if $I_p$ is not $A$-graded, there is
some degree $b \in B$ with all monomials of degree $b$ contained in
$I_p$.  Now homogeneous polynomials of degree $b$ are in one-to-one
correspondence with vectors in $k^{d_b}$.  The bijection takes the basis
vector with a one in the row corresponding to $x^u$ and zeroes
elsewhere to $x^u$, and is defined on other vectors by linear
extension.  Homogeneous polynomials of degree $b$ contained in $I_p$
are those corresponding to the image of the map $\sigma: k^{n_b}
\rightarrow k^{d_b}$ given by $\sigma: x \mapsto M_bx$.  Thus if all
monomials of degree $b$ are in $I_p$, $M_b$ must have full rank, which
means that there is a maximal minor which does not vanish.

While these equations for $H_A$ are not binomial, it follows from the
work of Peeva and Stillman \cite{PS2} that each irreducible component
of the scheme is given by binomial equations.  The work of Eisenbud
and Sturmfels on binomial ideals \cite{ES} now implies that the
radical of the ideal defining each component is also a binomial ideal,
and so the reduced structure on each irreducible component is a toric
variety.  We denote by $\tilde{H}_A$ the underlying reduced scheme of
$H_A$.

The main result of this section is:

\begin{theorem}\label{conniffconn}
The toric Hilbert scheme $H_A$ is connected if and only if the flip
graph of $A$ is connected.
\end{theorem}

The remainder of this section builds up to the  proof of Theorem
\ref{conniffconn}. By the support of a point $v \in \mathbb A^n$ we
mean $supp(v):=\{i : v_i \neq 0 \}$. In what follows we assume some
familiarity with toric varieties, such as that given in  \cite{Ewald}
or \cite{Fulton}. 
 
Corollary 2.6 of \cite{ES} says that every prime binomial ideal
determines a (not necessarily normal) toric variety.  The next lemma
gives a property of such varieties.  When $Q$ is a prime ideal of $S$
we denote by $V(Q)$ the the zero set of $Q$ in $\mathbb A^n$.

\begin{lemma} \label{torussupport}
Consider the point configuration $\{p_1,\dots,p_n\} \subseteq \mathbb
Z^d$ and its toric ideal $Q=ker(\phi:k[x_1,\ldots,x_n] \rightarrow
k[t^{p_1},\dots,t^{p_n}])$ which is a prime binomial ideal.  Let $v_1$
and $v_2$ be two points in $V(Q)\subseteq \mathbb A^n$.  Then $v_1$
and $v_2$ lie in the same torus orbit of $V(Q)$ if and only if they
have the same support.
\end{lemma}

\begin{proof}
  The dense torus in $V(Q)$ is $V(Q) \cap (k^*)^n$, and the action of
  this torus on $V(Q)$ is by coordinate-wise multiplication.  It thus
  follows that if $v_1$ and $v_2$ are in the same torus orbit, they
  have the same support.
  
  Suppose $v_1,v_2 \in V(Q)$ have the same support.  If this support
  is the entire set $\{1,\dots,n\}$, then define
  $u_i=(v_1)_i/(v_2)_i$.  Then if $x^a-x^b$ is a binomial in $Q$,
  $u^a-u^b=\frac{v_1^a}{v_2^a}-\frac{v_1^b}{v_2^b}=
  \frac{1}{v_2^av_2^b}(v_1^av_2^b-v_1^bv_2^a)=0$,   
  so $u$ is in $V(Q) \cap (k^*)^n$, and so $v_1$ and $v_2$ are in the
  same torus orbit.

Suppose now that $v_1$ and $v_2$ have the same support $\tau
\subsetneq \{1,\dots,n\}$. Since $v_1$ and $v_2$ are in $V(Q)$, this
means that there is no binomial in $Q$ of the form $x^a-x^b$ where
$supp(a) \subseteq \tau$ and $supp(b) \not \subseteq \tau$.  This is
because if such a binomial were in $Q$, we will have $v_i^b=0$ for
$i=1,2$, and $v_i^a \neq 0$ for $i=1,2$, which contradicts $v_1,v_2
\in V(Q)$.  This means that there is no affine dependency between
$\{p_i : i \in \tau \}$ and $\{p_i : i \not \in \tau \}$.  But this
implies that $conv(p_i : i \in \tau)$ is a face of $conv(p_i : 1 \leq
i \leq n)$, and if $p_j \in conv(p_i : i \in \tau)$, then $j \in
\tau$.  This means that $v_1$ and $v_2$ lie in an invariant toric
subvariety, and so by a similar argument to above are torus
isomorphic.
\end{proof}

The action of $(k^*)^n$ on $A$-graded ideals gives an action of
$(k^*)^n$ on $\tilde{H}_A$.  The $n$-torus acts by mapping $v\in
\tilde{H}_A$ to $t \cdot v$ via the map $(t \cdot v)_u=t^uv_u$.We will
refer to this action as the $n$-torus action.  There is also a torus
action on a point for every irreducible component of the reduced toric
Hilbert scheme it belongs to.  We will refer to these actions as the
ambient torus actions.  We note that these torus actions are usually
different from the $n$-torus action, as each of the finitely many
irreducible components of $H_A$ has only finitely many ambient torus
orbits, but there can be an infinite number of $n$-torus orbits.  An
example of this situation is given in Theorem 10.4 of \cite{GBCP}.
The $n$-torus orbit is, however, contained inside all ambient torus
orbits.

\begin{corollary} \label{toruscontainment}
Let $v$ be a point on $\tilde{H}_A$.  Then the $n$-torus orbit of $v$ is
contained in any ambient torus orbit of $v$.
\end{corollary}

\begin{proof}
It is straightforward to see that $t \cdot v$
lies in every irreducible component of $\tilde{H}_A$ in which $v$ does (this
follows from the fact that $S[l]/(lI_v+(1-l)(I_{t\cdot v}))$ is a flat
$k[l]$ module).  All points in the $n$-torus orbit of $v$ have the
same support, and thus lie on the same ambient torus orbit by Lemma
\ref{torussupport}.
\end{proof}

Fix an irreducible component $V$ of $\tilde{H}_A$.  Since $V$ is a
projective toric variety, there is a polytope $P$ corresponding to
$V$.  An ambient torus orbit of a point $v \in \tilde{H}_A$
corresponds to a face of $P$.  In the case of the coherent component,
this polytope is the state polytope of $I_A$.  Over the course of the
next three lemmas, we show that the edges of $P$ correspond exactly to
flips.

\begin{lemma} \label{verticesaremono}
Vertices of $P$ correspond exactly to the  monomial $A$-graded ideals in V.
\end{lemma}
\begin{proof}
Let $I$ be the ideal corresponding to a vertex $p$ of $P$.  The orbit
of $I$ under the ambient torus corresponding to $P$ is just the ideal
$I$.  By Corollary \ref{toruscontainment} the $n$-torus orbit of $I$
is contained in any ambient torus orbit, so $I$ is $n$-torus
fixed as well, and thus is a monomial ideal. 

For the other implication, let $I$ be a monomial $A$-graded ideal
corresponding to a point $v$ in $V$.  As a point in $X$, $v$ is
invariant under any scaling of its coordinates in any fashion, and
thus is invariant under any ambient torus action. It thus corresponds
to a vertex of $P$.
\end{proof}

\begin{lemma} \label{twoinismeansflip}
Let $I$ be an $A$-graded ideal.  If $I$ has exactly two initial
ideals, then $I$ is $n$-torus isomorphic to an ideal of the form
$J=\langle x^a-x^b, x^{c_1}, \ldots, x^{c_r} \rangle$.
\end{lemma}

\begin{proof}
Let $M_1$ and $M_2$ be the two initial ideals of $I$, and let
$\mathcal G$ be the universal Gr\"obner basis of $I$.  The set
$\mathcal G$ contains a reduced Gr\"obner basis for $I$ with respect
to a term order for which $M_1$ is the initial ideal, and so there
exist binomials $x^a-cx^b \in \mathcal G$ with $c \neq 0$ for which
$x^a$ is a minimal generator of $M_1$, $x^b \not \in M_1$.  Suppose
for all such binomials we have $x^a \in M_2$.  Then $M_1 \subseteq
M_2$ is an inclusion of distinct monomial $A$-graded ideals, which is
impossible.  So we conclude that there is some binomial
$x^{a_1}-c_1x^{b_1}\in \mathcal G$ with $c_1 \neq 0$, $x^{a_1} \in M_1
\setminus M_2$ and $x^{b_1} \in M_2 \setminus M_1$.

 Suppose there is some other binomial $x^{a_2}-c_2x^{b_2} \in \mathcal
G$ with $c_2 \neq 0$.  Without loss of generality we may assume that
$x^{a_2} \in M_1$ and $x^{b_2} \not \in M_1$.  We note that $(a_1-b_1)
\neq (a_2-b_2)$, as by Lemma \ref{UGB} the two binomials
$x^{a_1}-x^{b_1}$ and $x^{a_2}-x^{b_2}$ are Graver binomials, and they
must be distinct since $\mathcal G$ is the universal Gr\"obner basis
of $I$.  We can thus find a supporting hyperplane for
$pos(a_1-b_1,b_2-a_2)$, which intersects the cone only at the origin.
This implies the existence of a vector $w$ which satisfies $w \cdot
(a_1-b_1) >0$ and $w \cdot (b_2-a_2) > 0$.  Let $M=in_w(I)$.  Then
$x^{a_1} \in M$, and $x^{b_2} \in M$, so $M \neq M_1$, and $M \neq
M_2$.  This means that $I$ has a third initial ideal, which
contradicts our assumption, and so we conclude that
$x^{a_1}-c_1x^{b_1}$ is the only binomial in $\mathcal G$.

Pick $i\in supp(b_1)$.  Define $\lambda_i=\frac{1}{c_1}$, and
$\lambda_j=1$ for $j \neq i$.  Then $\lambda I$ is in the desired
form.
\end{proof}

\begin{theorem} \label{edgesareflips}
Let $M_1$ and $M_2$ be monomial $A$-graded ideals corresponding to
vertices $p_1$ and $p_2$ of $P$.  $M_1$ and $M_2$ are connected by a
single flip if and only if there is an edge $e$ of $P$ connecting
$p_1$ and $p_2$. 
\end{theorem}

\begin{proof}
Suppose $p_1$ and $p_2$ are connected by an edge $e$.  Let $I$ be the
ideal corresponding to a point $p$ in the relative interior of $e$.
By Corollary \ref{toruscontainment} the $n$-torus closure of $p$ is
contained in $e$.  Thus $I$ has at most two initial ideals.  If $I$
had only one initial ideal, it would be a monomial ideal and thus
corresponds to a vertex of $P$, by Lemma \ref{verticesaremono}.  We
thus conclude that $I$ has exactly two initial ideals, $M_1$ and
$M_2$, corresponding to $p_1$ and $p_2$ respectively.  Now by Lemma
\ref{twoinismeansflip} $I$ is $n$-torus isomorphic to $J= \langle
x^a-x^b, x^{c_1}, \ldots, x^{c_r} \rangle$, where $x^a \in M_1
\setminus M_2$ and $x^b \in M_2 \setminus M_1$.  Since $J$ is
$A$-graded, $x^a-x^b$ is a Graver binomial.  Because $J$ has
initial ideals $M_1$ and $M_2$, it is their wall ideal $W_{a-b}$, and
so $M_1$ and $M_2$ are connected by a flip over $x^a-x^b$.

Conversely, suppose $M_1$ and $M_2$ are connected by a single flip.
Then there is an ideal $W_{a-b}=\langle x^a-x^b, x^{c_1}, \ldots,
x^{c_r} \rangle$ which has as its two initial ideals $M_1$ and $M_2$.
Let $J$ be an $A$-graded ideal which is isomorphic to $W_{a-b}$ under
the ambient torus corresponding to $P$.  Let $x^d$ be a minimal
generator of $M_1$, with $d \neq a$, and $x^d-x^e$ the corresponding
Graver binomial with $x^e \not \in M_1$.  Then $x^d \in W_{a-b}$, and
thus $x^d \in J$, as the ambient torus action preserves the monomials
in an $A$-graded ideal.  So $J$ contains all minimal generators of
$M_1$ and $M_2$ except $x^a$ and $x^b$.  Suppose $J$ has a minimal
generator $x^{\alpha}-cx^{\beta}$, where $x^{\alpha}-x^{\beta}$ is a
Graver binomial, $x^{\alpha},x^{\beta} \not \in J$, and $\alpha, \beta
\neq a,b$.  Without loss of generality we may assume that $x^{\alpha}
\in M_1$. If $x^{\beta} \not \in M_1$ then $x^{\alpha} \in W_{a-b}$ by
the definition of $W_{a-b}$, and thus also $x^{\alpha} \in J$.  We
thus conclude that $x^{\beta} \in M_1$.  But this means there exist
$\alpha^{\prime} \leq \alpha$, $\beta^{\prime} \leq \beta$, such that
$x^{\alpha^{\prime}}$ and $x^{\beta^{\prime}}$ are minimal generators
of $M_1$.  Since $x^{\alpha}$ and $x^{\beta}$ have disjoint support,
we cannot have $\alpha^{\prime}=\beta^{\prime}=a$, so at least one of
$x^{\alpha^{\prime}}$ and $x^{\beta^{\prime}}$ is in $W_{a-b}$.  But
this means at least one of $x^{\alpha}$ and $x^{\beta}$ is in $J$,
giving a contradiction.  Hence the only binomial minimal generator of
$J$ is of the form $x^a-c^{\prime}x^b$, so as in the proof of Lemma
\ref{twoinismeansflip} $J$ is $n$-torus isomorphic to $W_{a-b}$.  We
thus see that all ambient torus closures of $W_{a-b}$ are the same as
the $n$-torus closure, and so $p_1$ and $p_2$ are connected by an
edge.
\end{proof}

\begin{proof}[Proof of Theorem \ref{conniffconn}]
It suffices to show that the reduced scheme $\tilde{H}_A$ is connected
if and only if the flip graph of $A$ is connected.  Since passing to
an initial ideal is a flat deformation, each irreducible component
contains a monomial $A$-graded ideal. It thus suffices to show that
all monomial $A$-graded ideals lie in the same connected component of
$\tilde{H}_A$ if and only if the flip graph is connected.  The ``if''
direction follows from the fact that if $I_1$ and $I_2$ are connected
by a single flip, then they are both initial ideals of a single wall
ideal $W_{a-b}$, and so lie in the same connected component of
$\tilde{H}_A$.  The ``only-if'' direction follows from Lemmas
\ref{verticesaremono} and \ref{edgesareflips}, which imply that the
flip graph restricted to an irreducible component of $\tilde{H}_A$ is
the edge skeleton of a polytope whose vertices are the monomial
$A$-graded ideals in that component, and so is connected.  As the
intersection of two irreducible components of $\tilde{H}_A$ contains a
monomial $A$-graded ideal by Gr\"obner deformation, this means that if
$\tilde{H}_A$ is connected, the flip graph of $A$ is connected.
\end{proof}
 
\section{Toric Hilbert Schemes of Arbitrarily High Dimension from 
Curves in $\mathbb P^4$} 

In this section we exhibit toric Hilbert schemes of arbitrarily high 
dimensions for which the associated toric varieties are curves in 
${\mathbb P^4}$. When $A$ has corank one, its Graver basis
consists of precisely one binomial $x^a-x^b$, and the flip graph of $A$
has only the two vertices $I= \langle x^a \rangle$ and $I' = \langle
x^b \rangle$ which are connected by the flip $x^a-x^b$. Hence $H_A$ is
one-dimensional and connected. All $A$-graded ideals of a corank two
matrix are coherent \cite{GP} which implies that the flip graph of $A$
is connected since it coincides with the edge graph of the state
polytope of $I_A$. In this case, $H_A$ has exactly one irreducible
component which is two dimensional and smooth \cite{PS1}. The toric
Hilbert scheme of a corank three matrix is at least three dimensional
since the irreducible component containing the coherent $A$-graded
ideals has dimension three. In contrast to the results in coranks one 
and two, Theorem~\ref{highdim} gives a family of two by five matrices
of corank three whose toric Hilbert schemes can have arbitrarily high  
dimensions. The projective toric variety of each matrix in the family
is a curve in ${\mathbb P^4}$. Note that both the corank $n-d$ and the
number of columns $n$ are fixed for these matrices.

\begin{theorem}\label{highdim}
For each $j \in {\mathbb N} \backslash \{0\}$, the toric Hilbert
scheme $H_{A(j)}$ of 
$$A(j) = \left   ( \begin{array}{ccccc} 1 & 1 & 1 & 1 & 1 \\ 0 & 1 &
    3+3j & 4+3j & 6+3j \end{array} \right )$$ 
has an irreducible component of dimension at least $j$.
\end{theorem}

These matrices were motivated by Example~5.11 in \cite{ST}, and 
the theorem was inspired by computer experiments on their flip graphs.
We first define the following monomial ideals and sets of binomials
that will be used in the proof of Theorem~\ref{highdim}.
For each $j \in {\mathbb N} \backslash \{0\}$, let \\

$\begin{array}{ll}
P_j = \langle c^2e, bc, a^2e, ace, ae^{j+2} \rangle,&
R_j = \langle a^5c^j, a^8c^{j-1}, \ldots, a^{5+3(j-1)}c \rangle,\\
Q_j = \langle be^{j+1}, a^2c^{j+1}, b^4e^j, c^{j+2} \rangle,&
S_j = \langle b^7e^{j-1}, b^{10}e^{j-2}, \ldots, b^{7+3(j-1)} \rangle
\end{array}$
\begin{center} and \end{center}
$\begin{array}{l}
{\mathcal P}_j = \{ c^2e-d^3, bc-ad, a^2e-b^2d, ace-bd^2,
ae^{j+2}-c^jd^3 \},\\
{\mathcal Q}_j = \{be^{j+1}-c^{j+1}d,
a^2c^{j+1}-b^3e^j, b^4e^j-a^3c^jd, c^{j+2}-ae^{j+1}\},\\
{\mathcal R}_j = \{a^{5+3t}c^{j-t} - b^{6+3t}e^{(j-1)-t}, \, \,
t=0,1,\ldots, j-1 \},\\
{\mathcal S}_j = \{ b^{7+3t}e^{(j-1)-t}-a^{6+3t}c^{(j-1)-t}d, 
\,\,t=0,1,\ldots, j-1 \}.
\end{array}$\\

\begin{lemma} \label{gbasis}
The ideal $M_j = P_j + Q_j + R_j + S_j$ is the initial ideal of
$I_{A(j)}$ with respect to the weight vector $w = (1,1,2,0,2)$. 
\end{lemma}

\begin{proof} By computing the $A(j)$-degree of both terms in each
binomial of ${\mathcal G}_j := {\mathcal P}_j \cup {\mathcal Q}_j \cup
{\mathcal R}_j \cup {\mathcal S}_j$, it can be seen that ${\mathcal
G}_j$ is a subset of $I_{A(j)}$. It can also be checked 
that for each binomial in
${\mathcal G}_j$, the positive term is the leading term with respect
to $w = (1,1,2,0,2)$. Hence $M_j = \langle in_{w}(g) : g \in
{\mathcal G}_j \rangle$ is contained in the initial ideal of
$I_{A(j)}$ with respect to $w$ and 
no generator of $M_j$ is redundant. The monomial ideal $M_j$ will
equal $in_{w}(I_{A(j)})$ if ${\mathcal G}_j$ is the reduced
Gr\"obner basis of $I_{A(j)}$ with respect to $w$. 
Consider the elimination order $x,y \succ a,b,c,d,e$ refined by 
the graded reverse lexicographic order $x > y$ on the first 
block of variables and the weight vector $w$ on the second 
block of variables. Then the reduced Gr\"obner basis of $I_{A(j)}$ 
with respect to $w$ is the intersection of the reduced Gr\"obner 
basis of \[ J(j) := \langle
a-x, b-xy, c-xy^{3+3j}, d-xy^{4+3j}, e-xy^{6+3j} \rangle \] with
respect to $\succ$ with $k[a,b,c,d,e]$ (see Algorithm 4.5 in
\cite{GBCP}). By a laborious check it can be shown that the reduced
Gr\"obner basis of $J(j)$ with respect to $\succ$ is \\

$\begin{array}{l}
{\mathcal G}_j \cup \{x-a, ya-b, ybd-ae, yc-d, yd^2-ce, y^2d-e,\\
yb^{3t+2}e^{j-t}-a^{3t+1}c^{(j+1)-t},\, t=0,\ldots,j,\\
y^2b^{3t+1}e^{j-t}-a^{3t}c^{(j+1)-t}, \, t=0,\ldots,j, \\
y^{3l}b^{3t}e^{p_l-t} - a^{3t-1}c^{p_l-t+1},\,\,l=1,\ldots,j,
\,\,t=1,\ldots, p_l:=(j+1)-l,\\ 
y^{3l+1}b^{3t-1}e^{p_l-t} - a^{3t-2}c^{p_l-t+1},\,\,l=1,\ldots,j,\,\,
t=1,\ldots, p_l,\\  
y^{3l+2}b^{3t-2}e^{p_l-t} - a^{3t-3}c^{p_l-t+1},\,\,l=1,\ldots,j,\,\,
t=1,\ldots, p_l \}. 
\end{array}$\\
\end{proof}

\begin{lemma}\label{highval} 
For each $j \in {\mathbb N} \backslash \{0\}$ the 
monomial $A(j)$-graded ideal $M_j$ from Lemma~\ref{gbasis} has exactly
$2j+4$ flippable binomials.
\end{lemma}

\begin{proof}
We will show that the binomials in ${\mathcal Q_j} \cup
{\mathcal R}_j \cup {\mathcal S}_j$ are flippable 
for $M_j$ while those in ${\mathcal P_j}$ are not. 
In order to show that a binomial $x^a-x^b$ is flippable for
$M_j$ we need to show that every $S$-polynomial (monomial in our case)
formed from the binomial $x^a-x^b$ (with $x^a$ as leading term) 
and a minimal generator $x^c$ of
$M_j$ different from $x^a$ reduces to zero modulo $W_{a-b} = \langle
x^c : c \neq a,\, x^c$ a minimal generator of $M_j \rangle + 
\langle x^a-x^b \rangle$. 

We first consider ${\mathcal R_j}$.
A binomial $a^{5+3t}c^{j-t} - b^{6+3t}e^{(j-1)-t}$ in ${\mathcal R}_j$
can form a non-trivial $S$-pair ($S$-monomial) with (i) $c^2e$,
(ii) $bc$, (iii) $a^2e$, (iv) $ace$, (v) $ae^{j+2}$, (vi) $a^2c^{j+1}$,
(vii) $c^{j+2}$ and (viii) a monomial $a^{5+3l}c^{j-l}$ from $R_j$
such that $t \neq l$. The remaining generators of $M_j$ (except 
$a^{5+3t}c^{j-t}$ itself) are relatively prime to $a^{5+3t}c^{j-t}$
and so the $S$-pairs formed reduce to zero by Buchberger's first 
criterion. We consider each case separately.

(i) The $S$-monomials formed from $c^2e$ and $a^{5+3t}c^{j-t}-
b^{6+3t}e^{(j-1)-t}$ are $b^{6+3t}c^pe^{j-t}$, 
$0 \leq t \leq j-1$, where $p=1$ if $j-t = 1$ and $p=0$ if $j-t > 
1$.\\
\indent \indent (a) If $t=0$, $b^6c^pe^j$ is a multiple of
  $b^4e^j \in Q_j$.\\
\indent \indent (b) If $t > 0$, $b^{6+3t}c^pe^{j-t}$
reduces to zero modulo $b^{7+3(t-1)}e^{j-t} \in S_j$.

(ii) The $S$-monomials formed from $bc$ are
$b^{7+3t}e^{(j-1)-t}$, $0 \leq t \leq j-1$ all of which 
lie in $S_j$ and hence reduce to zero modulo $W_{a-b}$.

(iii) The $S$-monomials between $a^2e$ and $a^{5+3t}c^{j-t}-
b^{6+3t}e^{(j-1)-t}$ are $b^{6+3t}e^{j-t}$ for $0 \leq t \leq j-1$.
If $t=0$, $b^6e^j$ is a multiple of $b^4e^j \in Q_j$, and
if $t > 0$ then $b^{6+3t}e^{j-t}$ is divisible by $b^{7+3(t-1)}e^{j-t}
\in S_j$.  

(iv) The $S$-monomials from $ace$ are $b^{6+3t}e^{j-t}$ for 
$0 \leq t \leq j-1$, all of which reduce to zero as in (iii).

(v) The monomial $ae^{j+2}$ gives 
$b^{6+3t}e^{2j+1-t}$ for $0 \leq t \leq j-1$, all of which reduce 
to zero modulo $be^{j+1} \in Q_j$.

(vi) From $a^2c^{j+1}$ we get 
$b^{6+3t}c^{t+1}e^{(j-1)-t}$, $0 \leq t \leq j-1$, all of which are 
multiples of $bc \in P_j$.

(vii) The $S$-monomials from $c^{j+2}$ are 
$b^{6+3t}c^{t+2}e^{(j-1)-t}$ which are also multiples of 
$bc \in P_j$ for $0 \leq t \leq j-1$.

(viii) For this last case, suppose first that $l < t \in
\{0,1,2,\ldots,j-1\}$. Then 
$lcm(a^{5+3l}c^{j-l},a^{5+3t}c^{j-t}) = a^{5+3t}c^{j-l}$ and the
$S$-monomial between $a^{5+3l}c^{j-l}$ and $a^{5+3t}c^{j-t}-
b^{6+3t}e^{(j-1)-t}$ is $b^{6+3t}c^{t-l}e^{(j-1)-t}$ which
is a multiple of $bc \in P_j$. If $l > t$, then the $S$-monomial is 
$a^{3(l-t)}b^{6+3t}e^{(j-1)-t}$ which is divisible by $a^2e \in P_j$ 
since $t < l \leq j-1$ and hence $t < j-1$.

Similarly, one can check that the binomials in 
${\mathcal Q_j} \cup {\mathcal S}_j$ are all flippable for $M_j$, 
which shows that $M_j$ has at least $2j+4$ flippable binomials.
To finish the proof, we argue that no binomial in ${\mathcal P_j}$ is 
flippable for $M_j$. 

(i) The $S$-binomial between $c^2e-d^3 \in {\mathcal P_j}$ and $bc \in
P_j$ is $bd^3$ which is not divisible by any generator of $M_j$. 

(ii) The binomials $bc-ad,\,a^2e-b^2d$ and $ace-bd^2 \in {\mathcal P_j}$
form the $S$-binomials $ade^{j+1}, \, b^3de^j$ and $b^2d^2e^j$
respectively with $be^{j+1} \in Q_j$. None of them can be divided by a 
minimal generator of $M_j$. 

(iii) The $S$-binomial of $ae^{j+2}-c^jd^3 \in {\mathcal Q_j}$ and
$a^2e \in P_j$ is $ac^jd^3$ which does not lie in $M_j$.

Hence $M_j$ has exactly $2j+4$ flippable binomials.
\end{proof}

\begin{proof}[Proof of Theorem~\ref{highdim}] 
The same proof as in Lemma~\ref{highval} shows that the generators 
of ${\mathcal I}(\mu_0, \ldots, \mu_{j-1}) := P_j + Q_j
+ \langle a^{5+3t}c^{j-t}-\mu_tb^{6+3t}e^{j-1-t}\,,t=0,\ldots,j-1 \rangle 
+ S_j$ form a Gr\"obner basis with respect to $w = (1,1,2,0,2)$
with initial ideal $M_j$, for every choice of scalars
$\mu_0,\ldots,\mu_{j-1}$ from the underlying field $k$.
Lemma~\ref{highval} proved this claim for the case $\mu_i = 1$, for an
$0 \leq i \leq j-1$ and $\mu_j = 0$ for all $j \neq i$. Since $M_j$ is
$A(j)$-graded, the $A(j)$-homogeneous ideal ${\mathcal
  I}(\mu_0,\ldots,\mu_{j-1})$ is also $A(j)$-graded for every choice
of scalars 
$\mu_0,\ldots,\mu_{j-1}$. Hence there is an injective polynomial map
from ${\mathbb A_k^j} \rightarrow H_{A(j)}$, such that $(\mu_0,
\ldots, \mu_{j-1})$ maps to the point on $H_{A(j)}$ corresponding
(uniquely) to ${\mathcal I}(\mu_0, \ldots, \mu_{j-1})$. Since
${\mathbb A_k^j}$ is irreducible, the image of this map lies entirely
in one irreducible component of the toric Hilbert scheme $H_{A(j)}$
and the dimension of this component is at least $dim({\mathbb A_k^j})
= j$.   
\end{proof}

\begin{remark}
In \cite{ST} it was conjectured that the maximum valency of a vertex 
in the state polytope of $I_A$ is bounded above by a function in just
the corank of $A$. As a particular case, it was also conjectured that 
if $A$ is of corank three, then every vertex in the state polytope of 
$I_A$ has at most four neighbors. This latter conjecture was recently 
disproved by Ho\c{s}ten and Maclagan \cite{HuT} who have found vertices 
with up to six neighbors. Lemma~\ref{highval} shows that even in
corank three, a vertex in the flip graph of $A$ can have arbitrarily
many neighbors.
\end{remark}

\section{Connection to the Baues Problem}

In this section we elaborate a connection between $A$-graded ideals
and the Baues problem for triangulations.  A good reference for all
forms of the Baues problem is \cite{reiner}.

A triangulation of a point configuration $\mathcal A=
\{a_1,\ldots,a_n\}\subseteq \mathbb R^d$ is a geometric simplicial complex
covering $\text{conv}(a_1,\ldots,a_n)$ with the vertices of each
simplex being a subset of $\mathcal A$.  Each simplex $\sigma$ is indexed by
the set $\{i: a_i \text{ is a vertex of }\sigma\}$.

\begin{figure} 
\epsfig{file=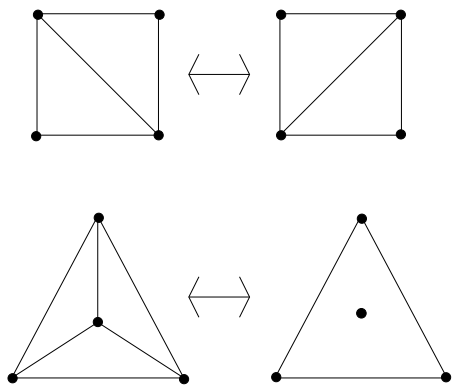, height=5cm}
\caption{Bistellar flips for triangulations of points in the plane}
\label{flipeg}

\end{figure}

A basic operation on triangulations of a point configuration is the
{\em bistellar flip}.  The two basic types of non-degenerate bistellar
flips in the plane are shown in Figure \ref{flipeg}.

Intuitively, a bistellar flip should be thought of as gluing in a
higher dimensional simplex, and then turning that simplex over and
viewing it from the other side.  This can be seen most clearly in the
second example in Figure \ref{flipeg}, which can be viewed as the
top and bottom of a tetrahedron.  The first example can also be
thought of as two opposite views of a tetrahedron.

More formally, a bistellar flip interchanges the two different
triangulations of a {\em circuit} (minimal affine dependence) of
$\mathcal A$.  Let $t$ be a circuit of the configuration $\mathcal A$,
and $T=\{i : t_i \neq 0\}$ be its support.  We denote by $T^+$ the set
$\{ i : t_i >0 \}$ and by $T^-$ the set $\{ i : t_i <0 \}$.  There are
exactly two triangulations of $C=conv(a_i : i \in T)$.  The first,
$C^+$, has $|T^+|$ simplices, which are the simplices indexed by the
sets $\{ T \setminus \{i\} : i \in T^+\}$.  The second, $C^-$, has
$|T^-|$ simplices, which are the sets in $\{ T \setminus \{i\} : i \in
T^-\}$.  The unique minimal non-face of $C^+$ ($C^-$) is $T^+$ ($T^-$).  If
$C$ is $d$-dimensional, and one of $C^+$ and $C^-$ is a subcomplex of
the triangulation $\Delta$, then a bistellar flip over the circuit $t$
involves replacing the subcomplex $C^+$ by $C^-$ or vice versa.

If $C$ is lower dimensional, we impose an additional condition for $t$
to be flippable.  By the {\em link} of a simplex $\sigma$ in a
simplicial complex $\Delta$ we mean the collection of simplices $\{
\tau : \tau \cap \sigma=\emptyset, \tau \cup \sigma \in \Delta \}$.
We say $t$ is flippable if $C^+$ (or $C^-$) is a subcomplex of
$\Delta$, {\em and} the link in $\Delta$ of every maximal simplex of
$C^+$ (respectively $C^-$) is  the same subcomplex $L$.  This second
condition is trivially satisfied if $C$ is $d$-dimensional, as the
link of every maximal simplex is the empty set.  A bistellar flip over
the circuit $t$ from $C^+$ to $C^-$ then involves replacing the
simplices $\{ l \cup \sigma : \sigma \in C^+, l\in L \}$ by the
simplices $\{ l \cup \tau : \tau \in C^-,l \in L\}$.

\begin{figure} 
\epsfig{file=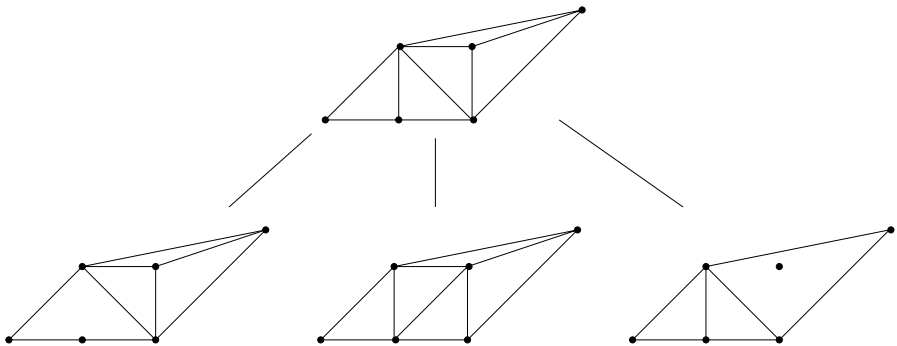, height=4.5cm}
\caption{Triangulations differing by bistellar flips}
\label{flipseg2}
\end{figure}

Examples of bistellar flips are shown in Figure \ref{flipseg2}.

We can form a graph, called the {\em Baues graph}, on the set of all
triangulations of a point configuration, with an edge connecting two
triangulations when they differ by a bistellar flip.  Figure
\ref{flipseg2} is a subgraph of the Baues graph for a particular
collection of six points in the plane.  An obvious question to ask is
whether the Baues graph is connected.  Santos recently answered this
question negatively \cite{Santos}, constructing a configuration of 324
points in $\mathbb R^6$ which has a disconnected Baues graph.

The rest of this section will relate the Baues graph to the flip graph
and the toric Hilbert scheme.  The connection is through the
following lemma, which is a special case of Theorem 10.10 in
\cite{GBCP}.  It links monomial $A$-graded ideals and triangulations
of $\mathcal A$, where $A$ is the matrix whose columns are the points
in $\mathcal A$, with an additional row of ones added.  We will denote
both the $i$th row of $A$ and the $i$th point of $\mathcal A$ by
$a_i$.  We adopt the notational convention that if $\sigma \subseteq
\{1,\dots,n\}$ is a set, then $x^{\sigma}=\prod_{i \in \sigma} x_i$.
The {\em Stanley-Reisner ideal} (see \cite{Stanley}) $I(\Delta)$ of a
simplicial complex $\Delta$ is the ideal generated by the monomials 
$x^{\sigma}$ where the sets $\sigma$ are the minimal
non-faces of $\Delta$.  Similarly, every squarefree monomial ideal $I$
in $S$ defines a unique simplicial complex $\Delta(I)$ on
$\{1,\dots,n\}$.

\begin{lemma} \label{triangs} \cite[Theorem~10.10]{GBCP}
Let $I$ be a monomial $A$-graded ideal.  Then $\Delta(rad(I))$, the
simplicial complex associated to $rad(I)$ via the Stanley-Reisner
correspondence, is a triangulation of $\mathcal A$.\hfill $\Box$
\end{lemma}

We can now state the main theorem of this section. 

\begin{theorem} \label{flipmeansbistellar}
Let $I$ be a monomial $A$-graded ideal and $x^a-x^b$ a flippable
binomial for $I$.  Then either $\Delta(rad(I))=\Delta(rad(I_{flip}))$, or they
differ by a bistellar flip.
\end{theorem}

The proof will be developed through the following series of lemmas.
We need to show that if $I_1$ and $I_2$ are monomial $A$-graded ideals
which differ by a flip, then either the radicals are the same, or
$\Delta(rad(I_1))$ and $\Delta(rad(I_2))$ differ by a bistellar flip.
This involves showing:

\begin{enumerate}
\item $t=a-b$ is a circuit of ${\mathcal A}$ (Lemma \ref{radsthesame}).
\item $C^+$ is a subcomplex of $\Delta(rad(I))$ with the link of all
maximal simplices of $C^+$ the same (Lemma \ref{ksubcomplex}).
\item $\Delta(rad(I_{flip}))$ differs from $\Delta(rad(I))$ exactly by
replacing $C^+$ and its link by $C^-$ and corresponding link.
\end{enumerate}

By a circuit of $A$ we mean a binomial $x^a-x^b$ such that $a-b$ is a
circuit of $\mathcal A$.  Note that all circuits are Graver binomials.

\begin{lemma}  \label{radsthesame}
Let $I$ be a monomial $A$-graded ideal, with $x^a-x^b$ a flippable binomial with $x^a
\in I$.  Then $x^b \in rad(I) \Leftrightarrow rad(I)=rad(I_{flip})$.
If $x^b \not \in rad(I)$, then $x^a-x^b$ is a circuit of $A$.
\end{lemma}
\begin{proof}
The implication $\Leftarrow$ is immediate in the first statement so we
need only show that $x^b \in rad(I)$ implies $rad(I)=rad(I_{flip})$.
Suppose $x^b \in rad(I)$.  Let $x^c$ be a minimal generator of
$I_{flip}$.  Then either $x^c$ is a minimal generator of $I$, $c=b$,
or $c=a+g$ for some $g$.  In each case $x^c \in rad(I)$, so
$rad(I_{flip}) \subseteq rad(I)$.  If the containment is proper, Lemma
\ref{triangs} gives a proper containment of Stanley-Reisner ideals of
triangulations of $\mathcal A$, which is not possible.  So we conclude
$rad(I)=rad(I_{flip})$.

For the second statement, suppose $x^a-x^b$ is not a circuit.  Then
there exists a circuit $x^c-x^d$ with $supp(c) \subseteq supp(a)$, and
$supp(d) \subseteq supp(b)$ where at least one of these inclusions is
proper.  Since $x^b \not \in rad(I)$, we must have $x^{supp(d)} \not
\in rad(I)$, and thus $x^d \not \in I$.  This implies $x^c \in I$, and
so, since we know $c \neq a$, $x^c \in I_{flip}$. This means
$x^{supp(c)} \in rad(I_{flip})$, and so $x^a \in rad(I_{flip})$.  But
this means, as above, that $rad(I_{flip})=rad(I)$, which in turn
implies that $x^b \in rad(I)$, contradicting the hypothesis.
\end{proof}

Let $I$ be a monomial $A$-graded ideal, with $x^a -x^b$ flippable,
where $x^a \in I$, $x^b \not \in rad(I)$.  Let $t=a-b$, and
$T=supp(t)$.  By Lemma \ref{radsthesame} we know that $t$ is a
circuit, so we can consider the triangulation $C^+=\{T \setminus \{i\} : i
\in T^+ \}$ of $C$.

\begin{lemma} \label{ksubcomplex}
Let $I$, $x^a-x^b$, $t$, and $C^+$ be as above.  Then $C^+$ is a
subcomplex of $\Delta=\Delta(rad(I))$, and there is a subcomplex of
$\Delta$ which is the common link of all maximal simplices of $C^+$.
\end{lemma}
\begin{proof}
$T^+$ is the only minimal non-face of $C^+$, so to show that $C^+$ is
 a subcomplex of $\Delta$, we need to show that $x^{T^+}$ is the only
 minimal generator of $rad(I)$ with support in $T$.

Suppose $x^c$ is a minimal generator of $rad(I)$, with $supp(c)
\subseteq T$.  Then there is some $l \geq 1$ such that $x^{lc} \in I$.
Write $c=a^{\prime}+b^{\prime}$, where $supp(a^{\prime}) \subseteq
supp(a)$, and $supp(b^{\prime}) \subseteq supp(b)$.  If
$supp(a^{\prime}) \neq supp(a)$, then $x^a$ does not divide
$x^{lc}$ and so $x^{lc}$ is in the wall ideal $W_{a-b}$.  We can
choose $\delta$ with $supp(\delta) \subseteq supp(b)$ so that $lc+
\delta=mb+a^{\prime}$ for some $m \geq 1$.  Since
$x^{mb+a^{\prime}}=x^{lc+\delta} \in W_{a-b}$, it follows that
$x^{ma+a^{\prime}} \in W_{a-b}$, because $x^a-x^b \in W_{a-b}$.  So
$x^{ ma+a^{\prime}} \in I_{flip}$, and there is thus some $p \geq m+1$
such that $x^{pa} \in I_{flip}$. This implies that $x^a \in
rad(I_{flip})$.  But, by Lemma \ref{radsthesame}, this means that
$rad(I_{flip})=rad(I)$, which in turn implies that $x^b \in rad(I)$,
contradicting our hypothesis.  So $supp(a^{\prime}) =
supp(a)=T^+$, and thus $x^{T^+} | x^c$.  This shows that $x^{T^+}$ is
the only minimal generator of $rad(I)$ with support in $T$, as
required.  From this we conclude that $C^+$ is a subcomplex of
$\Delta$.

We now show that every maximal simplex $\sigma \in C^+ \subseteq
\Delta$ has the same link.  We do this by showing that any simplex not
in the link of one maximal simplex of $C^+$ is not in the link of any
other maximal simplex of $C^+$.

 Suppose $\sigma \subseteq \{1,\dots,n\}$ is not a simplex in the link
 of a maximal simplex $\gamma$ of $C^+ \subseteq \Delta$, where
 $\gamma=T \setminus \{p\}$ for some $p \in T^+$ and $\sigma \cap T =
 \emptyset$.  Then $x^{\sigma \cup \gamma} \in rad(I)$, because
 $\sigma \cup \gamma$ is not a face of $\Delta$, and so there exists
 $l \geq 1$, and $x^{\alpha}$ a minimal generator of $I$ with $\alpha
 \neq a$, such that $x^{\alpha} | (x^{\sigma \cup \gamma})^l$.  Write
 $\alpha=a^{\prime}+b^{\prime}+\sigma^{\prime}$, where
 $supp(a^{\prime}) \subsetneq supp(a)$, $supp(b^{\prime}) \subseteq
 supp(b)$, and $supp(\sigma^{\prime}) \subseteq \sigma$.  Choose
 $\delta$ with $supp(\delta) \subseteq supp(a)$ such that
 $\alpha+\delta=ma+b^{\prime}+\sigma^{\prime}$ for some $m \geq 0$.
 Then because $x^{\alpha} \in W_{a-b}$, we have $x^{\alpha+\delta} \in
 W_{a-b}$, and so $x^{mb+b^{\prime}+\sigma^{\prime}}$ is in $W_{a-b}$
 and thus in $I$.  So
 $x^{supp(b)\cup supp(\sigma^{\prime})} \in rad(I)$.  Let $\tau$ be
 another maximal simplex of $C^+$, so $\tau=(\gamma \cup \{p\})
 \setminus \{p^{\prime}\}$ for some $p^{\prime} \in T^+$.  Then $supp(b)
 \cup supp(\sigma^{\prime}) \subseteq \tau \cup \sigma$, and so
 $x^{\tau \cup \sigma} \in rad(I)$, and thus $\sigma$ is not in the link
 of $\tau$ in $\Delta$.  This shows that every maximal simplex $\sigma
 \in C^+ \subseteq \Delta$ has the same link, as required.

\end{proof}

\begin{proof}[Proof of Theorem \ref{flipmeansbistellar}]
If $x^b \in rad(I)$ then $rad(I)=rad(I_{flip})$ by Lemma
\ref{radsthesame}, and so $\Delta(rad(I))=\Delta(rad(I_{flip}))$.

Suppose $x^b \not \in rad(I)$.  Then Lemma \ref{radsthesame} implies
that $t=a-b$ is a circuit of $\mathcal A$.  By Lemma \ref{ksubcomplex} $C^+$ is
a subcomplex of $\Delta(rad(I))$ with each maximal simplex of $C^+$
having the same link in $\Delta(rad(I))$.  It remains to show that
$\Delta(rad(I_{flip}))$ is the result of performing a bistellar flip
on $\Delta(rad(I))$.

Let $\Delta^{\prime}$ be
the result of performing the bistellar flip on $\Delta(rad(I))$ over
$t$, and let $M$ be the Stanley-Reisner ideal of $\Delta^{\prime}$.

We claim that M is the squarefree monomial ideal generated by
 $x^{supp(b)}$, all the generators of $rad(I)$ except $x^{supp(a)}$, and
 also all monomials of the form $x^{\sigma}$, such that $supp(a) \subseteq
 \sigma$, and $\sigma \setminus (T \cap \sigma)$ is not in the link of
 the maximal simplices of $C^+$.  Let $\alpha \subseteq
 \{1,\dots,n\}$.  Then $\alpha$ is a face of $\Delta^{\prime}$ exactly
 when either $\alpha$ is a face of $\Delta$ and $T^- \not \subseteq
 \alpha$, or $\alpha= T^+ \cup \tau \cup \gamma$, where where $\tau
\subsetneq
 T^-$, and $\gamma$ is in the link of the maximal simplices of $C^+$.
 This means that $\beta \subseteq \{1,\dots,n\}$ is not a face of
 $\Delta^{\prime}$ exactly when either $T^- \subseteq \beta$, or
 $\beta$ is not a face of $\Delta$ and $\beta \neq T^+ \cup \tau \cup
 \gamma$ for any $\tau \subsetneq T^-$ and $\gamma$ in the link of the
 maximal simplices of $C^+$. This proves the claim.

We now show that $rad(I_{flip}) \subseteq M$.  Let $x^{\alpha}$ be a
minimal generator of $I_{flip}$ such that $x^{supp(\alpha)}$ is a
minimal generator of $rad(I_{flip})$.  If $x^{\alpha}$ is also a
minimal generator of $I$, then $x^{supp(\alpha)}$ is in the square
free ideal generated by all the generators of $rad(I)$ except
$x^{supp(a)}$, so $x^{supp(\alpha)} \in M$.  Since $x^{supp(b)} \in
M$, the only case left to consider is $\alpha=a+g$ for some $g \neq 0$
with $b \not \leq g$.  Write $g=a^{\prime}+b^{\prime}+\gamma$, where
$supp(a^{\prime}) \subseteq supp(a)$, $supp(b^{\prime}) \subsetneq
supp(b)$, and $supp(\gamma) \cap T = \emptyset$.  Choose $\delta$ so
that $\delta+a^{\prime}=la+\tilde{a}$ for some $l \geq 0$, where
$supp(\tilde{a})=T^+\setminus \{p\}$ for some $p \in T^+$.  Since
$x^{\alpha}$ is a minimal generator of $I_{flip}$ different from
$x^a$, it is in $W_{a-b}$.  It thus follows that $x^{a+g+\delta} \in
W_{a-b}$, and so, because $x^a-x^b \in W_{a-b}$, we have 
$x^{(l+1)b+\tilde{a}+b^{\prime}+\gamma} \in W_{a-b}$ and thus in $I$.  
Since $supp((l+1)b+\tilde{a}+b^{\prime})=T \setminus \{p\}$,
$x^{(T\setminus \{p\}) \cup supp(\gamma)} \in rad(I)$ and thus
$supp(\gamma)$ is not in the link of the maximal simplices of $C^+$.
Because $supp(\gamma)=supp(\alpha) \setminus T$, this means
$x^{\alpha} \in M$, and therefore $rad(I_{flip}) \subseteq M$.

Now because $\Delta(rad(I_{flip}))$ and $\Delta^{\prime}$ are both
triangulations of $\mathcal A$, this inclusion cannot be proper. So
$M=rad(I_{flip})$, and thus $\Delta(rad(I_{flip}))$ is the result of
performing a bistellar flip on $\Delta(rad(I))$.
\end{proof}

\section{Toward a disconnected toric Hilbert scheme}

We conclude with some results that suggest the existence of a toric
Hilbert scheme. As mentioned earlier, Santos 
\cite{Santos} has recently constructed a six dimensional point
configuration with 324 points for which there is a triangulation that
admits no bistellar flips. Hence this configuration has a disconnected
Baues graph. By the results in \cite{berndpreprint} and the previous
section, every monomial $A$-graded ideal $I$ is supported on a  
triangulation of $\mathcal A$ via the correspondence $I \mapsto
\Delta(rad(I))$, and if two monomial
$A$-graded ideals are adjacent in the flip graph of $A$, then either
they are supported on the same triangulation or on two triangulations
that are adjacent in the Baues graph of $A$. Just as for monomial
$A$-graded ideals, there is a notion of coherence for triangulations
of ${\mathcal A}$. Every {\em coherent} triangulation of $\mathcal A$ 
(often called a {\em regular} triangulation in the literature) supports
at least one monomial $A$-graded ideal, and at least 
one of these ideals is coherent (see Chapter 8 in \cite{GBCP}). On the
other hand, Peeva has shown that if a triangulation of $\mathcal A$ is
{\em non-coherent}/{\em non-regular} then there may be no monomial
$A$-graded ideal whose radical is the Stanley-Reisner ideal of this 
triangulation (see Theorem~10.3 in \cite{GBCP}). Hence in order for 
Santos' example to lift to an example of a disconnected toric Hilbert 
scheme, it suffices to show that there is a monomial $A$-graded
ideal whose radical is the Stanley-Reisner ideal of his isolated
(non-regular) triangulation. 
A straightforward search for such a monomial $A$-graded ideal from
his $6 \times 324$ matrix is, however, a daunting computational
endeavor. Nonetheless, Santos' disconnected Baues graph seems to be
evidence in favor of a disconnected flip graph.

Recall that every coherent monomial $A$-graded ideal has at least
$n-d$ neighbors in the flip graph of $A$.
We say that a monomial $A$-graded ideal is {\em flip deficient} if its 
valency in the flip graph of $A$ is strictly less than $n-d$. All flip
deficient monomial $A$-graded ideals are necessarily
non-coherent. Before Santos constructed an isolated triangulation,
discrete geometers provided several examples of flip deficient
triangulations (triangulations with valency less that $n-d$ in the
Baues graph) as evidence in support of the existence of a disconnected
Baues graph. We provide examples of flip deficient monomial $A$-graded
ideals. 

\begin{theorem}
For each matrix $A(n) := [1 \,\,2\,\,3\,\,7\,\,8\,\,9\,\,a_7 \cdots
a_n]$ with $a_i \in {\mathbb N}$ and $9 < a_7 < \cdots < a_n$,  
there is a monomial $A(n)$-graded ideal with at most $n-3 < n-1 =
corank(A(n))$ flips.   
\end{theorem}

\begin{proof}
For the matrix $A = [1\,\,2\,\,3\,\,7\,\,8\,\,9]$, the monomial ideal 
$J = \langle x_1x_5, x_2x_4, x_1x_4, x_1x_2, x_4x_6, x_2x_6, x_1x_6,
x_3x_4, x_2^2x_3, x_1x_3, x_2x_5^2, x_2^2x_5, x_1^2,\\ x_3^2, x_2^4,
x_3x_5^3, x_4^2x_5^2, x_4^3, x_5^6, x_4x_5^4
\rangle$ is $A$-graded. The flippable binomials of $J$ are
$x_5^6-x_3x_6^5, \, x_2x_6-x_3x_5$ and 
$x_3^2-x_2^3$. In this example, there are 2910 monomial $A$-graded
ideals in total and the flip graph of $A$ is connected.

Consider the monomial ideal $J' = J + \langle x_7, \ldots, x_n \rangle
\subseteq k[x_1,\ldots,x_n]$ and a degree $b \in {\mathbb N}A(n) =
{\mathbb N}A = {\mathbb N}$. All the monomials in $k[x_1,\ldots, x_n]$
of $A(n)$-degree $b$ that are divisible by at least one of $x_7,
\ldots, x_n$ are in $J'$ by construction. Among the monomials in
$k[x_1, \ldots, x_6]$ of degree $b$ (there is at least one such since
$b \in {\mathbb N}A$), there is precisely one that is not in $J$ and
hence not in $J'$ and hence $J'$ is $A(n)$-graded.  If $x^a - x^b \in
k[x_1, \ldots, x_6]$ is flippable for $J'$ then $in_{x^a \succ x^b}(
\langle x^a - x^b \rangle + \langle x^c : x^c$ minimal generator of
$J$, $c \neq a \rangle + \langle x_7, \ldots, x_n \rangle)$ =
$J'$. The only non-trivial $S$-pairs that are produced during this
calculation are those between $x^a-x^b$ and a monomial minimal
generator $x^c$ of $J$. Since the resulting initial ideal equals $J'$,
it follows that $in_{x^a \succ x^b} ( \langle x^a - x^b \rangle +
\langle x^c : x^c$ minimal generator of $J$, $c \neq a \rangle) = J$
and hence $x^a - x^b$ is flippable for $J$. So $x^a - x^b$ must be one
of the three flippable binomials of $J$.  Additionally, each of the
minimal generators $x_7, \ldots, x_n$ of $J'$ provides a flippable
binomial and hence $J'$ has $3+(n-6) = n-3$ flippable binomials.
\end{proof}

\begin{remark}
We have not found matrices of corank three with flip deficiency in our
experiments. However, flip deficiency occurs in corank four.
Consider 
$A = [3\,\,6\,\,8\,\,10\,\,15]$ and its monomial $A$-graded ideal
\[ \langle ae, bd, ab^2, be, a^2, d^2, e^2, b^3, abc^2 \rangle. \]
The neighboring monomial $A$-graded ideals are: \\
$\langle ae, bd, ab^2, be, de, a^2, d^2, e^2, b^3 \rangle$ from $de-abc^2$,\\
$\langle ae, bd, ab^2, b^2e, a^2, d^2, e^2, b^3, acd, abc^2 \rangle$
from $acd-be$, and \\
$\langle ae, bd, ab^2, ad^2, be, b^2c, a^2, d^3, e^2, b^3, abc^2, d^2e
\rangle$ from $b^2c-d^2$.
\end{remark}

\begin{remark}
  The above computations were made using two different programs.
  Starting with a monomial initial ideal of the toric ideal $I_A$ one
  can compute all monomial $A$-graded ideals in the same connected
  component as this initial ideal by using the results in Section 2 to
  calculate all the neighbors of a monomial $A$-graded ideal.  This
  computation can be done using the program TiGERS \cite{HuT} with 
  the command {\tt tigers -iAe filename} where {\tt filename} is the
  standard input file for {\tt TiGERS} with the data of the matrix
  $A$.  In order to find all monomial $A$-graded ideals, we resort to
  a second program (available from the authors) that first computes
  the Graver basis of $A$ and then systematically constructs weakly
  $A$-graded monomial ideals by choosing one monomial from each Graver
  binomial to be in the ideal (cf. Lemma~\ref{UGB}). The program then
  compares the Hilbert series of each such ideal against that of an
  initial ideal of $I_A$ to decide if it is A-graded. Comparing the
  total number of ideals produced by the two programs gives a
  convenient way to decide if the flip graph is connected.
\end{remark}

We conclude with an algorithmic issue concerning the enumeration of
all $A$-graded monomial ideals in the same connected component as a
fixed one. The main program in TiGERS enumerates the vertices of the
state polytope of $I_A$ by using the {\em reverse search} strategy of
Avis and Fukuda \cite{AF}, which requires only the current vertex to
be stored at any given time.  The input to the program is any one
monomial initial ideal of $I_A$ from which the program reconstructs
all the others without needing to consult the list of ideals it has
already found.  An essential requirement of this algorithm is a method
by which the input ideal can be distinguished from any other monomial
initial ideal of $I_A$ by considering only the edges of the state
polytope. This is done in TiGERS as follows:

Suppose $M_1$ and $M_2$ are two monomial initial 
ideals of $I_A$ induced by the weight vectors $w_1$ and $w_2$
respectively. Let ${\mathcal G}_1$ and ${\mathcal G}_2$ be the
corresponding reduced Gr\"obner bases of $I_A$. Then for each facet
binomial $x^a-x^b$ in ${\mathcal G}_1$ we have $w_1 \cdot (a-b) >
0$ and for each facet binomial $x^{\alpha}-x^{\beta} \in {\mathcal
  G}_2$ we have $w_2 \cdot (\alpha-\beta) > 0$. The reduced Gr\"obner
bases ${\mathcal G}_1$ and ${\mathcal G}_2$ coincide if and only if
each facet binomial $x^{\alpha}-x^{\beta}$ of ${\mathcal G}_2$
satisfies the inequality $w_1 \cdot (\alpha - \beta) > 0$.  Suppose
the input is a fixed initial ideal of $I_A$.  By the previous
observation, every other monomial initial ideal of $I_A$ will have a
mismarked facet binomial with respect to this term order and hence can
be distinguished from the input ideal. The following example shows
that monomial $A$-graded ideals cannot be distinguished by checking
the orientation of their flippable binomials.

\begin{example}
Consider $A=[3\,\,4\,\,5\,\,13\,\,14]$ and its non-coherent monomial 
$A$-graded ideal
\[ M = \langle cd^5, c^2e^3, be, d^9, b^2, c^3, a^6, bd, ae^2, ad^3,
ac^2, a^2d, a^2b, bc, a^3e, a^3c \rangle. \]
The flippable binomials of $M$ are $ae^2-cd^2$, $c^3-a^5$ and
$d^9-ce^8$. With respect to the weight vector $w =
(0,0,1,20,22)$, each of these flippable binomials has its positive 
term as leading term and hence $M$ cannot be distinguished from 
$in_{w}(I_A)$ by checking whether its flippable binomials are
mismarked with respect to $w$.
\end{example}

\section{Acknowledgments}

We would like to thank Bernd Sturmfels for helpful conversations.

\end{document}